\documentclass[11pt]{article}
\PassOptionsToPackage{obeyspaces}{url}
\usepackage[hidelinks]{hyperref}
\usepackage{bookmark}
\bookmarksetup{open,numbered,depth=5}
\usepackage{amsfonts}
\usepackage{amsmath}
\usepackage{amssymb, amsthm, enumitem}
\usepackage{breakurl}
\usepackage{pgf,tikz}
\usepackage{bm}
\usetikzlibrary{calc}

\usepackage{etoolbox} 
\patchcmd{\thebibliography}{\leftmargin\labelwidth}{\leftmargin\labelwidth\addtolength\itemsep{-0.1\baselineskip}}{}{}

\usepackage{mathtools}
\usepackage{xstring}

\newcommand*\samethanks[1][\value{footnote}]{\footnotemark[#1]}

\author{Ting-Wei Chao\thanks{Department of Mathematical Sciences, Carnegie Mellon University, Pittsburgh, PA 15213, USA\@. Supported in part by U.S.\ taxpayers through NSF CAREER grant DMS-1555149.} \and Zichao Dong\samethanks}

\title{A simple proof of the Gan--Loh--Sudakov conjecture}
\date{}

\usepackage[nameinlink]{cleveref}

\newtheorem{theorem}{Theorem}
\newtheorem{lemma}[theorem]{Lemma}

\newcommand*{\eqdef}{\stackrel{\mbox{\normalfont\tiny def}}{=}} % definition by equality                                      
\newcommand*{\N}{\mathbb{N}}                                    % Natural numbers
\newcommand*{\Z}{\mathbb{Z}}                                    % Integers
\newcommand*{\Pg}{\Phi_{\text{good}}}
\newcommand*{\Pb}{\Phi_{\text{bad}}}
\newcommand*{\Og}{\Omega_{\text{good}}}

\crefname{enumi}{step}{steps}
\crefname{part}{part}{parts}

\allowdisplaybreaks[4]

\begin{document}

\maketitle

\begin{abstract}
	We give a new unified proof that any simple graph on $n$ vertices with maximum degree at most $\Delta$ has no more than $a\binom{\Delta+1}{t}+\binom{b}{t}$ cliques of size $t \ (t \ge 3)$, where $n = a(\Delta+1)+b \ (0 \le b \le \Delta)$. 
\end{abstract}

\section{Introduction}

For a positive integer $t \ge 3$, let $k_t(G)$ be the number of cliques of size $t$ in a simple graph $G = G(V, E)$. In \cite{gan_loh_sudakov}, Gan, Loh, and Sudakov asked how large $k_t(G)$ can be for graphs with maximum degree at most $\Delta$. They made a conjecture, which we henceforth refer to as the \emph{GLS Conjecture}, that $k_t(G)$ is maximized by a disjoint union of $a$ cliques of size $\Delta+1$ and one clique of size $b$, where $|V| = a(\Delta+1) + b$ for $0 \le b \le \Delta$. Moreover, they proved in \cite{gan_loh_sudakov} that 
\[
\text{the GLS Conjecture holds for $t = 3$ \quad $\Longrightarrow$ \quad the GLS Conjecture holds for $t \ge 4$. }
\]
The proof is an application of the Lov\'{a}sz version of the famed Kruskal--Katona theorem (see \cite{frankl}). 

Later on, Chase proved that the GLS Conjecture holds for $t = 3$ in \cite{chase}, and hence resolved the GLS Conjecture completely. In this short note we present a new proof of the GLS conjecture that works for all $t \ge 3$ uniformly without using the Kruskal--Katona theorem. The proof can be viewed as a simplification and a generalization of Chase's proof in \cite{chase}. We prove the following statement: 

\begin{theorem} \label{thm:GLS}
	Let $G$ be a simple graph on $n$ vertices with maximum degree at most $\Delta$. For any integer $t \ge 3$, if $n = a(\Delta+1)+b$ where $a, b \in \Z$ and $0 \le b \le \Delta$, then $k_t(G) \le a\binom{\Delta+1}{t}+\binom{b}{t}$. 
\end{theorem}

For every simple graph $G = G(V, E)$, write $u \sim v$ if $uv$ is an edge, and $u \nsim v$ if $uv$ is a nonedge. We denote by $\overline{N(v)} \eqdef \{v\} \cup \{u \in V : u \sim v\}$ the closed neighborhood of $v$. Let $T_v$ be the set of all $t$-cliques intersecting $\overline{N(v)}$. The proof of \Cref{thm:GLS} relies on the following lemma:  

\begin{lemma} \label{lem:neighb_removal}
	For any integer $t \ge 3$, if $G = G(V, E)$ is a simple graph, then 
	\[
	\sum_{v \in V} |T_v| \le \sum_{v \in V} \binom{\deg(v)+1}{t}. 
	\]
\end{lemma}

This note is organized as follows: We first show that \Cref{thm:GLS} follows from \Cref{lem:neighb_removal}, and then prove \Cref{lem:neighb_removal} in a separate section. 

\begin{proof}[Proof of \Cref{thm:GLS} assuming \Cref{lem:neighb_removal}] 
	Fix $t \ge 3$ and $\Delta \in \N_+$, and let $G$ be an $n$-vertex graph. Then there exists $v \in V$ such that $|T_v| \le \binom{\deg(v)+1}{t}$, by \Cref{lem:neighb_removal}. 
	
	We induct on $n$. The base case is obvious, as \Cref{thm:GLS} is trivially true for $n = 0, 1, \dotsc, \Delta+1$. Suppose \Cref{thm:GLS} is true for $n-1, n-2, \dotsc, n-\Delta-1$. Then we have that
	\[
	k_t(G) \le \begin{cases}
		\binom{\deg(v)+1}{t} + a\binom{\Delta+1}{t} + \binom{b-\deg(v)-1}{t}, \quad &\text{when $b \ge \deg(v)+1$}, \\
		\binom{\deg(v)+1}{t} + (a-1)\binom{\Delta+1}{t} + \binom{b+\Delta-\deg(v)}{t}, \quad &\text{when $b < \deg(v)+1 \le b+\Delta+1$}.  
	\end{cases}
	\]
	Since the sequence $\left\{\binom{n}{t}\right\}_{n \ge 0}$ is convex, we have that $\binom{\deg(v)+1}{t} + \binom{b-\deg(v)-1}{t} \le \binom{b}{t}$ when $b \ge \deg(v)+1$, and $\binom{\deg(v)+1}{t}+\binom{b+\Delta-\deg(v)}{t} \le \binom{\Delta+1}{t}+\binom{b}{t}$ otherwise. We conclude that $k_t(G) \le a\binom{\Delta+1}{t}+\binom{b}{t}$. 
\end{proof}

\section{Proof of \texorpdfstring{\Cref{lem:neighb_removal}}{Lemma 2}}

Define the set 
\[
\Phi \eqdef \{(u, x_1, \dotsc, x_t) \in V^{t+1} : \text{$x_1, \dotsc, x_t$ form a $t$-clique in $G$, and $u \sim x_i$ for some $i \in [t]$}\}. 
\]
%Suppose $(v,  x_1, \dotsc, x_t) \in \Phi$, then at least one of $x_1, \dotsc, x_t$ lies in $\overline{N(v)}$
Observe that each $(v, x_1, \dotsc, x_t) \in \Phi$ consists of a vertex $v \in V$ and a $t$-clique $x_1 \dotsb x_t \in T_v$. Since for every $t$-clique in $G$, there are $t!$ ways to label its $t$ vertices as $x_1, \dotsc, x_t$, we have that 
\begin{equation} \label{eq:T_Phi}
	|\Phi| = t! \sum_{v \in V} |T_v|. 
\end{equation}

For each tuple $(u, x_1, \dotsc, x_t) \in \Phi$, the vertices $u, x_1, \dotsc, x_t$ are not necessarily distinct. However, there are at least $t$ distinct vertices among $u, x_1, \dotsc, x_t$, because $x_1, \dotsc, x_t$ form a $t$-clique. For every tuple $(u, x_1, \dotsc, x_t) \in V^{t+1}$, we call it \emph{good} if $u, x_1, \dotsc, x_t$ are distinct, and \emph{bad} otherwise. Let
\begin{align*}
	\Pg &\eqdef \{(u, x_1, \dotsc, x_t) \in \Phi : \text{$(u, x_1, \dotsc, x_t)$ is good}\}, \\
	\Pb &\eqdef \{(u, x_1, \dotsc, x_t) \in \Phi : \text{$(u, x_1, \dotsc, x_t)$ is bad}\}. 
\end{align*}
Then $\Pg$ and $\Pb$ partition $\Phi$. 

Fix $v \in V$. If $(v, x_1, \dotsc, x_t) \in \Pb$, then $v, x_1, \dotsc, x_t$ are vertices of a $t$-clique in $G$, where exactly one $x_i$ happens to be $v$. There are $t$ choices for this $x_i$, and at most $\binom{\deg(v)}{t-1}$ choices for the rest of the vertices $x_1, \dotsc, x_{i-1}, x_{i+1}, \dotsc, x_t$, and $(t-1)!$ choices for their possible permutations. Hence, 
\begin{equation} \label{eq:phi_good}
	|\Pb| \le \sum_{v \in V} t \cdot \binom{\deg(v)}{t-1} \cdot (t-1)! = t! \sum_{v \in V} \binom{\deg(v)}{t-1}. 
\end{equation}

To upper bound $|\Pg|$, we need to introduce the auxiliary set
\begin{align*}
	\Og \eqdef \{(w, y_1, \dotsc, y_t) \in V^{t+1} : \ &\text{$(w, y_1, \dotsc, y_t)$ is good, $w \sim y_i$ for all $i \in [t]$, } \\
	&\text{and $y_1, \dotsc, y_t$ contain a $(t-1)$-clique in $G$}\}. 
\end{align*}
For any fixed $v \in V$, if $(v, y_1, \dotsc, y_t) \in \Og$, then $y_1, \dotsc, y_t$ are distinct neighbors of $v$, and so
\begin{equation} \label{eq:omega_good}
	|\Og| \le t! \sum_{v \in V} \binom{\deg(v)}{t}. 
\end{equation}

We claim that
\begin{equation} \label{eq:phi_omega}
	|\Pg| \le |\Og|. 
\end{equation}
Assume that \eqref{eq:phi_omega} is established. From the combination of \eqref{eq:T_Phi}, \eqref{eq:phi_good}, \eqref{eq:omega_good}, and \eqref{eq:phi_omega}, we obtain
\begin{align*}
	t! \sum_{v \in V} |T_v| &= |\Phi| = |\Pb| + |\Pg| \le |\Pb| + |\Og| \\
	&\le t! \sum_{v \in V} \Bigg( \binom{\deg(v)}{t-1} + \binom{\deg(v)}{t} \Bigg) \\
	&= t! \sum_{v \in V} \binom{\deg(v)+1}{t}, 
\end{align*}
which concludes the proof of \Cref{lem:neighb_removal}. \qed

\begin{proof}[Proof of estimate \eqref{eq:phi_omega}]
When $\bm{u} \eqdef (u, x_1, \dotsc, x_t) \in \Pg$ or $\bm{w} \eqdef (w, y_1, \dotsc, y_t) \in \Og$, the induced subgraph $G[\bm{u}]$ or $G[\bm{w}]$ is connected and contains a $t$-clique. Consider any induced $(t+1)$-vertex subgraph $H$ of $G$ that is connected and contains a $t$-clique. Let $z_1, \dotsc, z_t$ be the vertices of the $t$-clique (choose arbitrary ones if there are several). Let $z^*$ be the remaining vertex of $H$. Assume without loss of generality that $z^* \sim z_1, \dotsc, z^* \sim z_k$, and $z^* \nsim z_{k+1}, \dotsc, z^* \nsim z_t$. Note that $t \ge 3$, we count for different values of $k$ the contribution of $H$ to $|\Pg|$ and $|\Og|$, respectively: 
\begin{itemize}
	\item $1 \le k \le t-2$. If $(u, x_1, \dotsc, x_t) \in \Pg$, then $u = z^*$ since the degree of $z^*$ in $H$ is less than $t-1$, and hence $\{x_1, \dotsc, x_t\} = \{z_1, \dotsc, z_t\}$. If $(w, y_1, \dotsc, y_t) \in \Og$, then $w \in \{z_1, \dotsc, z_k\}$, and hence $\{y_1, \dotsc, y_t\} = \{z^*, z_1, \dotsc, z_t\} \setminus \{w\}$. Such an $H$ contributes $t!$ and $k \cdot t!$ elements to $\Pg$ and $\Og$, respectively. 
	\item $k = t-1$. If $(u, x_1, \dotsc, x_t) \in \Pg$, then $\{x_1, \dotsc, x_t\} \supset \{z_1, \dotsc, z_{t-1}\}$, and hence $u \in \{z_t, z^*\}$. If $(w, y_1, \dotsc, y_t) \in \Og$, then $w \in \{z_1, \dotsc, z_{t-1}\}$, and hence $\{y_1, \dotsc, y_t\} = \{z^*, z_1, \dotsc, z_{t-1}\} \setminus \{w\}$. Such an $H$ contributes $2 \cdot t!$ and $(t-1) \cdot t!$ elements to $\Pg$ and $\Og$, respectively. 
	\item $k = t$. Then $H = K_{t+1}$. Such an $H$ contributes $(t+1)!$ elements to both $\Pg$ and $\Og$. 
\end{itemize}
The claimed estimate \eqref{eq:phi_omega} follows from the cases above. 
\end{proof}

\section*{Acknowlegements}

We are grateful to Boris Bukh for helpful discussions and suggestions during the preparation of this note. We thank an anonymous referee for valuable feedback on the earlier version of this note.

\bibliographystyle{plain}
\bibliography{proof_GLS}

\end{document}